# Fourier series-based algorithm for control optimization in pendulum capsule drive: an integrated computational and experimental study

Sandra Zarychta[a][1], Marek Balcerzak[a], Katarzyna Wojdalska[b], Rafał Dolny[b], Jerzy Wojewoda[a]

[a] *Division of Dynamics, Lodz University of Technology, Stefanowskiego 1/15, Lodz, Poland*

[b] *Lodz, Poland*

**Abstract**: Pendulum-driven systems have emerged as a notable modification of vibro-impact mechanisms, replacing the conventional mass-on-spring oscillator with a pendulum. Such systems exhibit intricate behavior resulting from the interplay of directional dynamics, pendulum motion, and contact forces between the designed device and the underlying surface. This paper delves into the application of a Fourier series-based greedy algorithm for control optimization in pendulum capsule drives, which hold potential for diverse scenarios, including endoscopy capsule robots, pipeline inspection, and rescue operations in confined spaces. The emphasis is placed on experimental studies involving prototype development to validate the system's efficacy with previous computational simulations. Empirical findings closely align (<2% loss) with numerical investigations, showcasing the pendulum capsule drive's ability to achieve average speeds of 2.48 cm/s and 2.58 cm/s for three and six harmonics, respectively. These results are reinforced by high-quality signal-tracking accuracy, which demonstrates resilience against potential disturbances during motion. The authors envision the Fourier series-based control optimization method as a significant step towards ensuring enhanced locomotion performance in discontinuous systems, effectively handling the non-linearities arising from dry friction.

## 1   Introduction

Lower gastrointestinal (GI) diseases, including cancer, Crohn's disease, ulcerative colitis, diverticular disease, GI bleeding, and polyposis syndromes, afflict numerous patients annually, with approximately one in five patients facing inaccurate diagnoses or underestimations during clinical assessments [1, 2]. Over the past two decades, capsule endoscopy has emerged as a potent technique for detecting and diagnosing GI disorders [3, 4]. By leveraging the inherent benefits of pill-sized devices that are equipped with functional modules like miniaturized cameras, capsule endoscopy provides wireless operation, painless procedures, enhanced safety, and quick positioning. These advantages make it a promising alternative to traditional endoscopes [2, 5]. The emergence of active capsule robotics has further augmented the precision and capabilities of endoscopic procedures, facilitating tasks such as biopsy, drug administration, and cancer identification [6–8]. Nevertheless, realizing these functionalities necessitates precise control over capsule locomotion within the gastrointestinal tract, posing a significant challenge for engineers in the field of capsule robotics [6, 9, 10]. Further investigations into control optimization, propulsion, and locomotion dynamics are vital. For instance, potential hindrances impeding the locomotion of capsule endoscopes during cancer detection in the intestines, such as circular folds, tumors, and lesions, must be addressed [2]. Examples of propulsion mechanisms include legged capsules [11–13], spiral-based capsules [14], and vibration capsules [15] [16].

---

[1] corresponding author: sandra.zarychta@dokt.p.lodz.pl

Current methods for facilitating capsule locomotion often rely on externally located driving mechanisms, including appendages such as arms, fins, or propellers. However, these mechanisms carry the inherent risk of causing trauma to the intestinal environment [1]. Inspired by the locomotion of inchworms [2], self-propelled mobile mechanisms driven by autogenous internal force and environmental resistance have garnered significant attention from applied mathematicians, experimentalists, and engineers. Their theoretical challenges as piecewise-smooth dynamical systems and broad applications in robotics make them particularly intriguing [17]. The significant upside of capsules that move solely due to gravity and friction with the surface has led to the emergence of an actively propelled capsule design with a smooth surface, which is both patient-friendly and clinically viable [1].

Within the literature, discontinuous capsule mechanisms have been investigated for their ability to navigate due to inertial forces and frictional environments. These can be grouped into three main categories based on their internal driving mechanisms: vibro-impact, pendulum-like, and vibration-driven systems [17]. Alternative categories based on internal driving mechanisms are also occasionally identified. Ito et al. [18] recently introduced impulse-driven capsules, presenting a compact, smooth-surfaced capsule powered by inertia and friction forces, suitable for medical procedures like inspection and drug delivery.

Initially, the more commonly studied vibro-impact systems merit discussion [19]. The vibro-impact mechanism is a self-propelled device driven by internal harmonic excitation, facilitating rectilinear motion while overcoming environmental resistance. Various research pursuits have explored comparative investigations [19] [20], modeling approaches [21] [22], vibration [23], [24], and multistability control techniques [25], frictional environment considerations [1] [2] [26], experimental validations [27] [28] [17], prototype development [17] [29], progression rates and optimization of energy consumption [30] [31], dynamic response continuations [21] [32], and near-grazing dynamics [33] [34]. For example, Zhang et al. observed that the vibro-impact method outperformed constant pulling [19]. Guo et al. conducted experimental investigations involving mesoscale vibro-impacting capsules, assessing their feasibility across diverse frictional environments. Through optimization of a standard-sized capsule prototype, they achieved a maximum average forward velocity of 8.49 mm/s, suggesting the potential for real-time and controllable examinations [28]. Several substantial works in the area were conducted by Liu et al. In paper [21], they determined that optimizing for maximum progression does not necessarily result in the most energy-efficient parameters, highlighting the need for a trade-off between progression and energy consumption in future optimizations. In work [25], the team investigated the control of multistability within a vibro-impact capsule system by implementing a position feedback controller. The system demonstrated forward and backward motion capabilities by altering the initial conditions. Utilizing the COCO platform, robustness analysis identified parameter regions conducive to effective control strategies. Moreover, they validated the obtained results on a physical prototype, revealing mainly periodic behavior and potential performance improvements through control parameter adjustments [27]. Last but not least, the team examined a capsule system propelled by a harmonic force using different friction models, optimizing force parameters for directed motion and speed enhancement [23].

Secondly, it is imperative to delve into examples of pendulum-like capsule drives. Through the substitution of the conventional mass-on-spring oscillator with a pendulum, a noteworthy modification of the vibro-impact capsule can be realized. In this context, the propulsion of the capsule results from the intricate interplay among directional dynamics, inertial forces deriving from the pendulum's swinging motion, and the contact force between the capsule and the underlying surface [35]. This configuration appears to introduce a level of complexity to the system dynamics, given the dependence of the contact force on the oscillations of the pendulum. The principles of periodic locomotion and the nonlinear dynamics inherent in a pendulum-driven capsule system have been investigated [35–37] and tested via simulations considering the control for the pendulum capsule system in Liu et al. (sinusoidal) [37–39] and our (Fourier) previous studies [40] [41].

Liu et al. investigated varying stiffness coefficients, identifying a region conducive to periodic motions. Incorporating selected parameters enhanced performance, increasing travel distance by 15.76% and reducing energy consumption by 16.46% [37]. In our previous research [40], it was discovered that using the Fourier-based algorithm led to over a 140% improvement in distance traveled compared to the optimized sinusoidal control. Despite maintaining a consistent signal amplitude, the Fourier-based method proved more efficient in optimizing capsule drives. The existing literature on pendulum-like driving mechanisms in robotic capsules is noticeably limited and predominantly theoretical. Addressing this gap constitutes a key objective of our ongoing research pursuits.

Regarding vibro-driven mechanisms, Nunuparov et al. examined the motion of a capsule-type mobile robot traversing a straight path on a rough horizontal surface using analytical and experimental methodologies [42]. Liu et al. studied the dynamic analysis and stick-slip effect of a capsule system, focusing on the influence of elasticity and viscosity. Their research revealed that the optimal selection of these coefficients enables desired forward motion and discusses stick-slip motion regions where nonlinearities dominate [43]. In their study, Du Nguyen et al. introduced a novel locomotion module for a capsule robot, employing a moving mass linked by a spring and stimulated by electromagnetic force to induce oscillatory motion [44]. La et al. conducted an experimental analysis on a vibration-driven locomotion system designed for capsule robots, uncovering the substantial impact of friction force on both speed and directional movement [45]. These results offer valuable insights to guide future design and operational considerations for similar systems. Liu et al. introduced an energy-preserving design to enhance efficiency in vibro-driven robotic (VDR) systems through a spring-augmented pendulum. Validation via comparative simulations and experimental verification underscores its effectiveness, positioning it as a pioneering prototype in VDR system research [46].

The potential applications of this technology in confined spaces, such as endoscopy, disaster rescue, or pipeline inspection, are explored, albeit requiring further refinement and study [37] [47]. This work focuses on developing a discontinuous, non-smooth pendulum capsule drive utilizing control optimization through an open-loop Fourier series-based algorithm started in [40]. Building upon the theoretical groundwork laid by us in the previous study, in which the Fourier-based control optimization method was assessed, this research applies this approach to optimize control in non-smooth mechanical systems, addressing discontinuities resulting from friction or impacts in the pendulum capsule drive [40].

In this paper, a modified implementation of the Fourier series-based algorithm realizing a greedy strategy is introduced. In the optimization procedure, the locally optimal solution for shorter time intervals is utilized as an initial condition for optimizing control over a longer time span, ultimately achieving the control function across the entire time domain as the final objective. This approach ensures the non-decreasing quality of solutions in subsequent repetitions of the procedure, which is confirmed by the numerical optimization of chosen numbers of harmonics. In order to validate the aforementioned findings, a prototype of the system and a test stand were constructed. The experimental validation of the proposed approach, including locomotion performance and signal-tracking accuracy, confirmed the effectiveness of the Fourier series-based method for nonlinear systems. This offers an alternative control strategy to existing methods. This study introduces a unique experimental approach that demonstrates the effectiveness of the Fourier series in addressing discontinuous systems, with potential implications for similar systems across various applications.

This manuscript is structured as follows: Section 2 provides a detailed description of the Fourier-series-based algorithm and its application of a greedy strategy for optimization. Section 3 introduces the pendulum capsule drive, including the mathematical model, device design, control strategy, parameter identification, numerical conditions, prototype and experimental setup, and measurement technique. Section 4 presents the results and subsequent discussion. Finally, the work is summarized, and conclusions are drawn.

## 2 Numerical control optimization – methodology

The methodology used to search for the optimal control in this research is built upon a modified, greedy version of the Fourier series-based algorithm. This chapter provides an overview of the approach, starting with a concise description of the original Fourier series-based method. It then proceeds to introduce the greedy variant of the algorithm, explaining the reasons for its creation and presenting how it operates.

### 2.1 The original Fourier series-based algorithm

To keep this paper self-contained, a brief introduction to the original Fourier series-based algorithm is provided. The description provided aims to provide enough information to comprehend the remaining content of the paper. However, readers seeking a more extensive and comprehensive description are encouraged to read the original work [40]. Several features of the method, which are only briefly mentioned here, have been thoroughly discussed and proven (where necessary) in the original paper.

Suppose that an optimal control function $u^*: [0, T] \rightarrow [u_{min}, u_{max}]$ is to be determined. Let $\tilde{u}$ be its Fourier approximation, as follows (1):

$$\tilde{u}(t) = \frac{a_0}{2} + \sum_{k=1}^{K} [a_k \cos(k\omega t) + b_k \sin(k\omega t)], \quad t \in [0, T] \quad (1)$$

where $a_0, a_1, \ldots, a_K, b_1, \ldots, b_K$ are Fourier coefficients [48] that are real constants, $\omega = 2\pi/T$ is the fundamental frequency, and $K$ is the number of harmonics in the signal. If $u^*$ satisfies the Dirichlet conditions and $K \rightarrow \infty$, then the Fourier series (1) converges to $u^*$ at all points where $u^*$ is continuous. Therefore, for assumed values of parameters $K$ and $\omega = 2\pi/T$, the approximation of the optimal control $u^*$ can be reduced to searching for the optimal coefficients $a_0, a_1, \ldots, a_K, b_1, \ldots, b_K$. In this manner, optimization of the function is reduced to a nonlinear programming problem.

However, the fundamental difficulty in such an approach is that the control must belong to the admissible set: $\forall_{t \in [0,T]} u(t) \in [u_{min}, u_{max}]$, and it is not known how to specify boundaries for the optimization of $a_0, a_1, \ldots, a_K, b_1, \ldots, b_K$ to ensure that this condition is fulfilled. To solve this problem, the following approach has been proposed. The formula (1) can be expressed in vector form (2):

$$\tilde{u}(t) = \frac{a_0}{2} + [a_1, b_1, \ldots, a_K, b_K][\cos(\omega t), \sin(\omega t), \ldots, \cos(K\omega t), \sin(K\omega t)]^T =$$
$$= \frac{a_0}{2} + \boldsymbol{H}[\cos(\omega t), \sin(\omega t), \ldots, \cos(K\omega t), \sin(K\omega t)]^T \quad (2)$$

where $\boldsymbol{H} = [a_1, b_1, \ldots, a_K, b_K]$ is a real vector in $\mathbb{R}^{2K}$. Since $\omega = 2\pi/T$ and $K$ are assumed a priori, $\boldsymbol{H}$ and $a_0$ are sufficient to provide a complete description of $\tilde{u}$. Out of these two quantities, $a_0/2$ is simply an additive constant: if the graph of $\tilde{u}(t)$ is taken into account, any change in $a_0$ only causes a shift of this graph along the y-axis. All other properties of $\tilde{u}$ are defined by $\mathbf{H}$.

Different features of the function $\tilde{u}$ are connected with the norm (length) of the vector $\mathbf{H}$, and others with its direction. To differentiate between both groups, it seems convenient to introduce two terms: *shape* and *span* of $\tilde{u}$. These two notions are illustrated in Fig. 1 and explained in the following paragraphs.

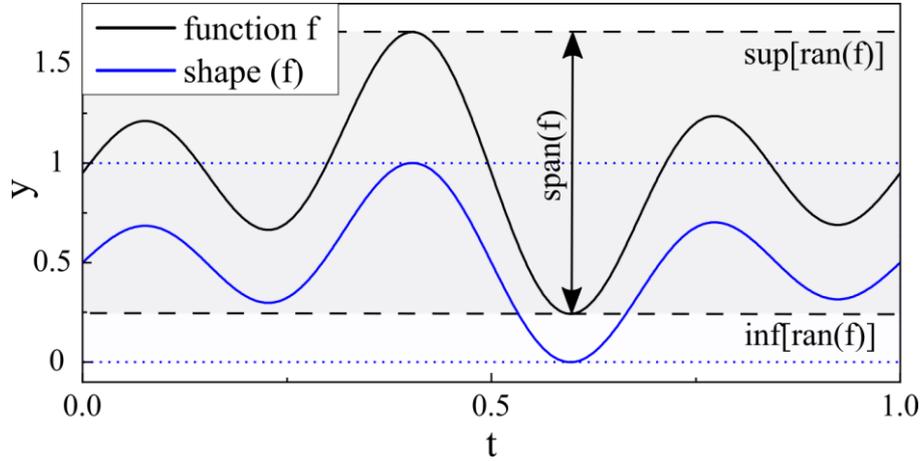

**Fig. 1** *Illustration of the notions of the shape and the span of a function f*

Let $D_{\tilde{u}} = [0, T]$ be the domain of $\tilde{u}$ and $ran(\tilde{u})$ denote its range (set of values). Then $\inf ran(\tilde{u})$ and $\sup ran(\tilde{u})$ are infimum and supremum of the range, respectively. The *span* of $\tilde{u}$ is defined as the following difference.

$$\mathrm{span}(\tilde{u}) = \sup ran(\tilde{u}) - \inf ran(\tilde{u}) \tag{3}$$

The definition (3) can be understood as the length of the smallest interval that is a superset of $ran(\tilde{u})$ – see **Fig. 1**. It can be noticed that multiplying $\tilde{u}$ by a positive real number causes the span to be multiplied by the same value, as formula (4) indicates. On the other hand, adding any constant to the function $\tilde{u}$ does not affect its span, as also depicted by formula (4).

$$\forall_{\substack{p>0 \\ q \in \mathbb{R}}} \mathrm{span}(p \cdot \tilde{u} + q) = p \cdot \mathrm{span}(\tilde{u}) \tag{4}$$

However, properties of $\tilde{u}$ which are not altered when $\tilde{u}$ is multiplied by a positive constant exist. Such multiplication, which graphically can be associated with "stretching" the graph of $\tilde{u}$ along the y-axis, leaves many of its features unchanged. Intervals of monotonicity, or argument values corresponding to maxima, minima, and points of inflection, as well as other similar properties of $\tilde{u}$, remain invariant under multiplication of $\tilde{u}$ by a positive number. The definition of *shape* presented in formula (5):

$$\mathrm{shape}(\tilde{u}) = \frac{\tilde{u} - \inf ran(\tilde{u})}{\mathrm{span}(\tilde{u})} \tag{5}$$

is intended to "move" and "stretch" the graph of $\tilde{u}$ so that it "fits" into the interval [0, 1] – see Fig. 1. Distinctly, formula (5) only makes sense if $\tilde{u}$ is not a constant function, otherwise the denominator equals 0. An obvious implication of equation (5) is that the shape of $\tilde{u}$ is invariant with respect to multiplication of $\tilde{u}$ by a positive constant or addition of any number to it.

$$\forall_{\substack{p>0 \\ q \in \mathbb{R}}} \mathrm{shape}(p \cdot \tilde{u} + q) = \mathrm{shape}(\tilde{u}) \tag{6}$$

Another consequence of formula (5) is that to have the function $\tilde{u}$ well-defined, it is sufficient to know its shape, span and infimum of its range.

$$\tilde{u} = \text{span}(\tilde{u}) \cdot \text{shape}(\tilde{u}) + \inf ran(\tilde{u}) \tag{7}$$

It is crucial to observe the connection between formulas (6) and (7), which depict important properties of the notions of shape and span, with equation (2) that expresses the function $\tilde{u}$ in terms of the constant vector $\boldsymbol{H}$ and the value $\frac{a_0}{2}$ (other terms, i.e., $\omega, K$, are assumed a priori). Formula (6) implies that multiplication of the vector $\boldsymbol{H}$ by a positive constant does not change the shape of the function $\tilde{u}$. Consequently, it is the direction of $\boldsymbol{H}$ that governs the shape of $\tilde{u}$ (along with parameters $\omega, K$).

Since $\boldsymbol{H} \in \mathbb{R}^{2K}$, its direction is well-defined by a point on a unit hypersphere of dimension $(2K - 1)$, which can be specified by $(2K - 1)$ spherical coordinates $\varphi_1, \ldots, \varphi_{2K-1}$ [49]. These spherical coordinates are bounded by the following intervals: $\varphi_1, \ldots, \varphi_{2K-2} \in [0, \pi]$ and $\varphi_{2K-1} \in [0, 2\pi)$.

$\widehat{\boldsymbol{H}} = [\hat{h}_1, \ldots, \hat{h}_{2K}] = \boldsymbol{H}/|\boldsymbol{H}|$ denotes the normalization of the vector $\boldsymbol{H}$. The values $\hat{h}_1, \ldots, \hat{h}_{2K}$ are related to the spherical components on the unit hypersphere by following equations (8).

$$\begin{aligned}
\hat{h}_1 &= \cos(\varphi_1) \\
\hat{h}_2 &= \sin(\varphi_1) \cos(\varphi_2) \\
&\ldots \\
\hat{h}_{2K-1} &= \sin(\varphi_1) \sin(\varphi_2) \ldots \sin(\varphi_{2K-2}) \cos(\varphi_{2K-1}) \\
\hat{h}_{2K} &= \sin(\varphi_1) \sin(\varphi_2) \ldots \sin(\varphi_{2K-2}) \sin(\varphi_{2K-1})
\end{aligned} \tag{8}$$

The function $\hat{u}$ is defined in formula (9).

$$\hat{u}(t) = \widehat{\boldsymbol{H}}[\cos(\omega t), \sin(\omega t), \ldots, \cos(K\omega t), \sin(K\omega t)]^T \tag{9}$$

It is evident that the function (9) can be obtained from (2) by setting $a_0 = 0$ and normalizing $\boldsymbol{H}$. Consequently, formula (6) shows that functions $\tilde{u}$ and $\hat{u}$ are of the same shape. The latter is well-defined by the values $K, \omega$ and the vector $\widehat{\boldsymbol{H}}$, whose components are directly related to the spherical coordinates $\varphi_1, \ldots, \varphi_{2K-1}$. This leads to the crucial conclusion: the shape of $\tilde{u}$ is well-defined by the angles $\varphi_1, \ldots, \varphi_{2K-1}$, which belong to known intervals. Hence, optimization of the shape of $\tilde{u}$ becomes a nonlinear programming problem involving optimization of the spherical coordinates $\varphi_1, \ldots, \varphi_{2K-1}$.

Although the essential step is finalized, formula (7) shows that it is not enough to know the shape of $\tilde{u}$ to have it well-defined: it is also necessary to specify its span and infimum of its range[2]. Since the assumed set of allowable controls is just an interval $[u_{min}, u_{max}]$ (see formula (1)), location of the range of $\tilde{u}$ in this set can be easily parametrized.[3] Parameters $p$ and $q$ can be defined according to formulas (10), (11) respectively.

$$p = \frac{\sup ran(\tilde{u}) - u_{min}}{u_{max} - u_{min}} \tag{10}$$

$$q = \frac{\sup ran(\tilde{u}) - \inf ran(\tilde{u})}{\sup ran(\tilde{u}) - u_{min}} = \frac{\text{span}(\tilde{u})}{\sup ran(\tilde{u}) - u_{min}} \tag{11}$$

Since $\tilde{u}$ is assumed to be non-constant and its range is a subset of $[u_{min}, u_{max}]$, the value $\sup ran(\tilde{u})$ satisfies the inequality $u_{max} \geq \sup ran(\tilde{u}) > u_{min}$, which according to formula (10) leads to $p \in$ (0, 1]. Analogously, since $\sup ran(\tilde{u}) > \inf ran(\tilde{u}) \geq u_{min}$, formula (11) yields $q \in (0, 1]$. Moreover, formulas (10), (11) can be transformed to the form (12), (13).

---

[2] In the context under consideration, infimum of the range of $\tilde{u}$ is simply equal to its minimum [48]. This results from the fact that $\tilde{u}$ is continuous (see formula (2)) and defined on the closed interval $[0, T]$.
[3] The presented method allows other types of control constraints too [40].

$$\sup ran(\tilde{u}) = p(u_{max} - u_{min}) + u_{min} = p \cdot u_{max} + (1 - p) \cdot u_{min} \quad (12)$$

$$\inf ran(\tilde{u}) = q[u_{min} - \sup ran(\tilde{u})] + \sup ran(\tilde{u}) = (u_{max} - u_{min})(1 - q)p + u_{min} \quad (13)$$

Equations (12), (13) prove that supremum and infimum of the range of $\tilde{u}$, and consequently its span, can be specified by the parameters $p, q \in (0, 1]$. This leads to the following, final conclusions concerning the optimization algorithm.

- Fourier approximation of the control function $\tilde{u}$ is well-defined by its shape, span and infimum of the range (see formula (7)).

- The shape of $\tilde{u}$ is well-defined in terms of spherical coordinates $\varphi_1, \ldots, \varphi_{(2K-1)}$ (formulas (8), (9)).

- The span of $\tilde{u}$ and the infimum of its range are well-defined by $p, q$ (formulas (12), (13)).

- The three points above imply that $\tilde{u}$ is well-defined by $(2K + 1)$ parameters: $\varphi_1, \ldots, \varphi_{2K-2} \in [0, \pi], \varphi_{2K-1} \in [0, 2\pi), p \in (0, 1], q \in (0, 1]$. Consequently, control optimization can be reduced to nonlinear programming problems of $(2K + 1)$ values.

Due to the fact that this optimization problem may involve multiple local extrema, it is solved using non-gradient methods, such as Differential Evolution (DE) [50] or Particle Swarm Optimization (PSO) [51].

## 2.2 Greedy variant of the Fourier series-based algorithm

A major drawback of evolutionary algorithms, such as Differential Evolution (DE), utilized in global optimization, is the challenge of solving problems with an increasing number of optimized variables. This issue can be observed in practical applications, as the number of iterations required to find an acceptable solution grows with the dimensionality of the underlying problem [52]. Obviously, the method presented in the previous section is vulnerable to such difficulties for large number of harmonics $K$ (see formula (2)), as the total number of parameters to be optimized equals $(2K + 1)$.

A possible way of overcoming issues connected to an increasing problem space in the DE algorithm is to split the problem into smaller parts and optimize each part separately [53]. In this subsection, a modification of the Fourier series-based algorithm that follows such a principle is presented. The approach applies a greedy strategy [54], in which locally optimal choices for shorter time intervals are used as a starting point for control optimization over a more extended period until the final target, i.e., control function over the whole-time domain $[0, T]$, is obtained. Such an idea is presented in **Fig. 2**. This section provides information about such methodology.

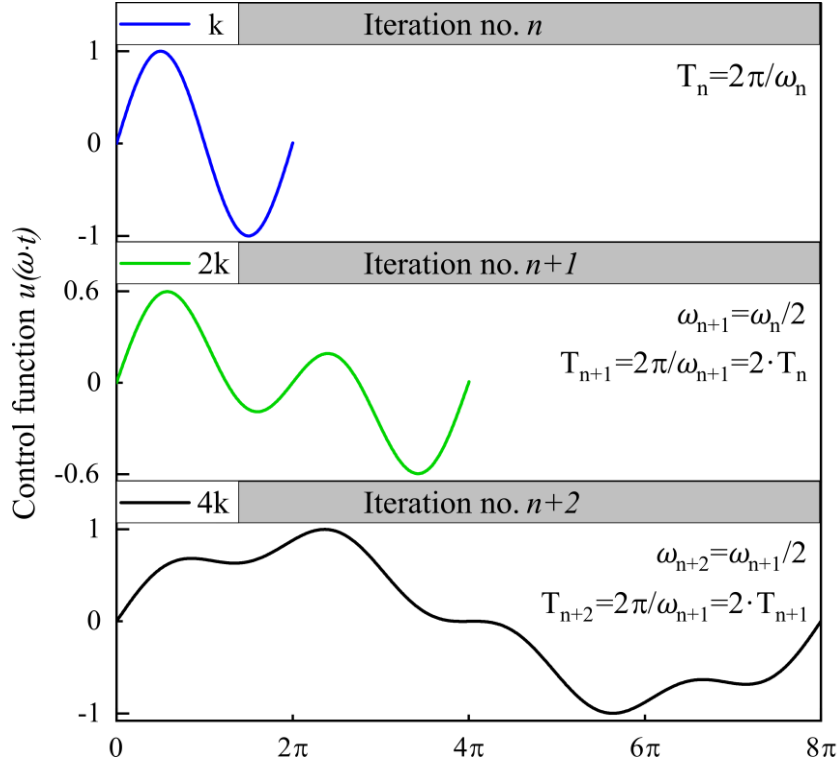

*Fig. 2 Illustration of the greedy variant of Fourier series-based control optimization algorithm*

Assuming that control optimization is going to be divided into $N$ steps, the length of the time interval over which the control is optimized will be doubled in each step. The time interval length in the $n$-th step is defined by the formula (14).

$$T_n = T \cdot 2^{n-p}, \quad n \in \{1, \dots, N\} \tag{14}$$

For instance, for $N = 3$ steps, $T_1 = T/4$, $T_2 = T/2$ and $T_3 = T$. Naturally, the first step is performed according to the description presented in the previous section.

For the step of a number $n$, the optimization begins over the interval $[0, T_n]$, and thus, the parameter $\omega$ (see formula (1)) attains the value $\omega_n = 2\pi/T_n$. If the control function expressed by equation (1) is to contain $k$ harmonics as it is optimized over $[0, T_n]$, their frequencies are $\omega_n, 2\omega_n, \dots, k\omega_n$ for a selected value $k$. Results of optimization are going to be the spherical coordinates $^n\varphi_1, \dots, ^n\varphi_{2k-1}$ (see formula (8)), and $^np$, $^nq$ (note formulas (10)-(13)), where the upper-left index means index of the step. These parameters enable calculation of the shape (equations (8), (9)) and span (formulas (12), (13)) of the control function over the interval $[0, T_n]$.

An important consideration is how to utilize the results obtained from the $n$-th step of control optimization over an interval $[0, T_n]$ to support estimation of the optimal control over a twice wider interval $[0, T_{n+1}] = [0, 2T_n]$? (Fig. 2). It is assumed that accuracy of optimal control estimation should be uniform, i.e., if the optimization takes place over an interval of a doubled length, the number of harmonics participating in the optimized control needs to be doubled as well. This raises the question about frequencies of subsequent harmonics. Since in the previous step of optimization the length of time interval was equal $T_n$, the fundamental frequency in the control function was $\omega_n = 2\pi/T_n$ and the control function contained harmonics of frequencies $\omega_n, 2\omega_n, \dots, k\omega_n$. In this particular scenario, the interval length doubles, and a new fundamental frequency emerges, amounting to half of the previously adopted value, as indicated by formula (15).

$$\omega_{n+1} = \frac{2\pi}{T_{n+1}} = \frac{2\pi/T_n}{2} = \frac{\omega_n}{2} \tag{15}$$

Frequencies of subsequent harmonics in the new iteration of the procedure are equal $\omega_{n+1}, 2\omega_{n+1}, \ldots, 2k\omega_{n+1}$ (following the assumption, the number of harmonics doubles when a twice longer time interval is considered). Formula (15) shows that the even harmonics from this list, i.e. $2\omega_{n+1}, 4\omega_{n+1}, \ldots, 2k\omega_{n+1}$ have already appeared in the control function in the previous step and their amplitudes (see formula (2)) have been optimized. However, odd harmonics $\omega_{n+1}, 3\omega_{n+1}, \ldots, (2k-1)\omega_{n+1}$ have not been present in the control function before. Therefore, to use the results of the previous step (control optimized over $[0, T_n]$ with even harmonics only $2\omega_{n+1}, 4\omega_{n+1}, \ldots, 2k\omega_{n+1}$) in the current iteration (optimization over the interval $[0, T_{n+1}] = [0, 2T_n]$ with all harmonics $\omega_{n+1}, 2\omega_{n+1}, \ldots, 2k\omega_{n+1}$), it is necessary to select such initial conditions for the optimization procedure that indicate the optimal result from the previous step. In this procedural framework, it is postulated that the control optimized in the present step is going to exhibit a performance that is at least comparable, if not superior, to its predecessor.

To practically realize this idea, it is necessary to analyze the parameters $\hat{h}_i$ and the corresponding values of spherical coordinates $\varphi_i$ that specify the shape of the optimized function (see formula (8)). From the previous step of the optimization procedure, with the index $n$, optimized values of ${}^n\varphi_1, \ldots, {}^n\varphi_{2k-1}$ and ${}^np, {}^nq$ are obtained. The former, i.e. spherical coordinates ${}^n\varphi_i, i \in \{1, \ldots, 2k-1\}$, specify shape of the optimized function by amplitudes of subsequent harmonics ${}^n\hat{h}_i, i \in \{1, \ldots, 2k\}$ – see formula (16), which depicts the optimized control function after the $n$-th step of optimization.

$$
\begin{aligned}
{}^n\hat{u}(t) &= {}^n\hat{h}_1 \cos(\omega_n t) + {}^n\hat{h}_2 \sin(\omega_n t) + {}^n\hat{h}_3 \cos(2\omega_n t) + \cdots + {}^n\hat{h}_{2k} \cos(k\omega_n t) = \\
&= {}^n\hat{h}_1 \cos(2\omega_{n+1} t) + {}^n\hat{h}_2 \sin(2\omega_{n+1} t) + {}^n\hat{h}_3 \cos(4\omega_{n+1} t) + \cdots + {}^n\hat{h}_{2k} \cos(2k\omega_{n+1} t)
\end{aligned} \quad (16)
$$

In the next step, the optimized function is described by the formula (17), in which the coefficients ${}^{n+1}\hat{h}_i, i \in \{1, \ldots, 4k\}$ are going to be optimized.

$$
{}^{n+1}\hat{u}(t) = {}^{n+1}\hat{h}_1 \cos(\omega_{n+1} t) + {}^{n+1}\hat{h}_2 \sin(\omega_{n+1} t) + {}^{n+1}\hat{h}_3 \cos(2\omega_{n+1} t) + \cdots + {}^n\hat{h}_{2k} \cos(k\omega_n t) \quad (17)
$$

Under the condition that equations (18) are met, the formula (17) can be employed to describe an equivalent function as (16). Therefore, formulas (18) describe the initial conditions of the next step of optimization which assure that the process starts from the optimized solution obtained in the previous step.

$$
\begin{gathered}
{}^{n+1}\hat{h}_1 = {}^{n+1}\hat{h}_2 = 0 \\
{}^{n+1}\hat{h}_3 = {}^n\hat{h}_1, \quad {}^{n+1}\hat{h}_4 = {}^n\hat{h}_2 \\
\ldots \\
{}^{n+1}\hat{h}_{4i} = {}^{n+1}\hat{h}_{4i+1} = 0 \\
{}^{n+1}\hat{h}_{4i+3} = {}^n\hat{h}_{2i+1}, \quad {}^{n+1}\hat{h}_{4i+4} = {}^n\hat{h}_{2i+2} \quad i \in \{1, \ldots, k\} \\
\ldots \\
{}^{n+1}\hat{h}_{4k-3} = {}^{n+1}\hat{h}_{4k-2} = 0 \\
{}^{n+1}\hat{h}_{4k-1} = {}^n\hat{h}_{2k-1}, \quad {}^{n+1}\hat{h}_{4k} = {}^n\hat{h}_{2k}
\end{gathered} \quad (18)
$$

Indisputably, the parameters ${}^{n+1}\hat{h}_i, i \in \{1, \ldots, 4k\}$ are not optimized directly, but through spherical coordinates. Comparison of equations (18) with the formula (8) leads to the following conditions for the angles.

$$
\begin{aligned}
&{}^{n+1}\varphi_1 = {}^{n+1}\varphi_2 = \frac{\pi}{2} \\
&{}^{n+1}\varphi_3 = {}^n\varphi_1, \quad {}^{n+1}\varphi_4 = {}^n\varphi_2 \\
&\quad \ldots \\
&{}^{n+1}\varphi_{4i} = {}^{n+1}\varphi_{4i+1} = \frac{\pi}{2} \\
&{}^{n+1}\varphi_{4i+3} = {}^n\varphi_{2i+1}, \quad {}^{n+1}\varphi_{4i+4} = {}^n\varphi_{2i+2} \\
&\quad \ldots \\
&{}^{n+1}\varphi_{4k-3} = {}^{n+1}\varphi_{4k-2} = \frac{\pi}{2} \\
&{}^{n+1}\varphi_{4k-1} = {}^n\varphi_{2k-1}
\end{aligned}
\qquad i \in \{1,\ldots,k\} \qquad (19)
$$

In summary, if the parameters ${}^n\varphi_1, \ldots, {}^n\varphi_{2k-2} \in [0, \pi]$, ${}^n\varphi_{2k-1} \in [0, 2\pi)$, ${}^n p \in (0, 1]$, ${}^n q \in (0, 1]$ are obtained at the $n$-th step of the optimization procedure, the next step should begin with the angles specified by formula (19) and the same values of values $p$ and $q$ as obtained in the previous step. This approach allows the optimization algorithm to automatically consider the previously optimized solution and ensures a non-decreasing quality of solutions throughout the optimization process.

Finally, it is essential to note that the formula (19) is applicable when the only constraint imposed on the control is the assumption $u(t) \in [u_{min}, u_{max}]$. However, in the later part of the paper, an additional constraint arises: $u(0) = 0$, which means that one of the amplitudes $\hat{h}_1$ is no longer independent and does not need to be optimized. In this scenario, assuming $\varphi_1 = \pi/2$ in all optimization steps, the formula (8) can be utilized to calculate all the amplitudes except for the first one, i.e. $\hat{h}_2, \ldots, \hat{h}_{2K}$ and to determine $\hat{h}_1$ from the conditions resulting from the additional constraints. Alternatively, modifying formulas (18) and (19) is an option, but it would lead to increased complexity and is not advisable.

## 3 Pendulum capsule drive

In the pendulum-like driven discontinuous capsule mechanisms, the swinging motion of the pendulum generates inertial forces. These forces, combined with the friction between the capsule and the surface beneath it, produce a horizontal force that moves the entire device in the desired direction. The system dynamics is evenly more intricate regarding the dependency of the contact load on the pendulum oscillations, which affects the friction force. **Fig. 3** illustrates the simplified conceptual motion of the pendulum capsule drive. A notable advantage of this system is that external moving parts, such as wheels, tracks, robotic arms, etc., are no longer required to produce the motion.

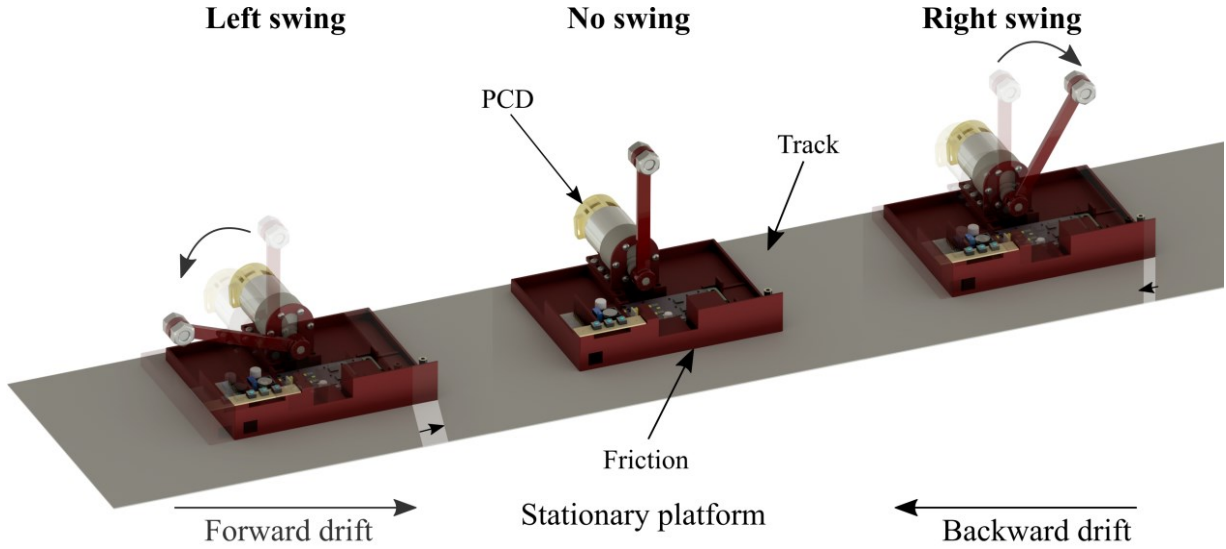

*Fig. 3 Simplified conceptual motion of the pendulum capsule drive (PCD)*

This section introduces the crucial equations governing the dynamics of the pendulum capsule drive, facilitating a deeper understanding of the interactions between contact load, inertial, and frictional forces. Following the theoretical background, the system's design is presented, along with a description of the electronic components and the system control approach. Details on the prototype, experimental setup, parameters identification, and numerical optimization and simulation conditions are then provided. Finally, the section concludes with an employed measurement technique to validate the system's behavior.

### 3.1 Mathematical model

This research focuses on the pendulum capsule drive, which scheme is illustrated in **Fig. 4**. The motion equations for this discontinuous system are thoroughly described in [40]. Thus, only a concise overview is provided here.

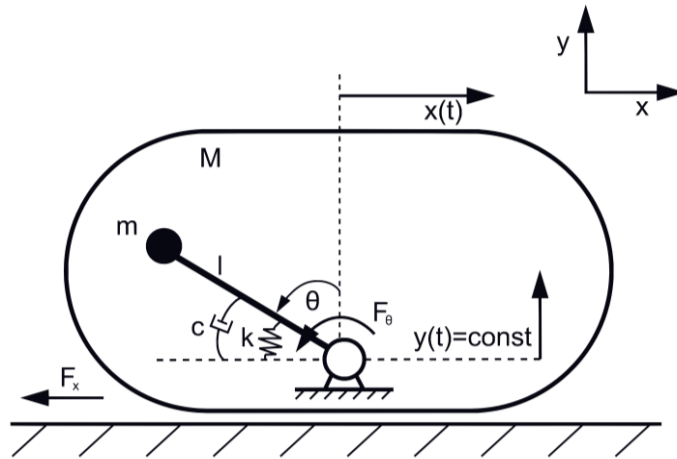

*Fig. 4 Scheme of the pendulum capsule drive. $M$ - mass of the capsule, $m$ – mass of the pendulum, $l$ - length of the pendulum, $\theta$ – pendulum angle, $k$ – spring stiffness, $c$ – damping coefficient, $F_\theta$ – external torque acting on the pendulum, $F_x$ – friction force, $x(t), y(t)$ – coordinates of the capsule*

The governing set of equations of the system in dimensional form are as follows:

$$(M + m)\ddot{x}(t) - ml\ddot{\theta}(t) \cos \theta(t) + ml\dot{\theta}^2(t) \sin \theta(t) = -F_x(t) \qquad (20)$$

$$ml^2\ddot{\theta}(t) - ml\ddot{x}(t)\cos\theta(t) = mgl\sin\theta(t) - k\theta(t) - c\dot{\theta}(t) + F_\theta(t) \tag{20'}$$

$$R_y(t) = \lambda_y(t) = (M+m)g - ml\ddot{\theta}(t)\sin\theta(t) - ml\dot{\theta}^2(t)\cos\theta(t) \tag{21}$$

$$R_x(t) = ml\ddot{\theta}(t)\cos\theta(t) - ml\dot{\theta}^2(t)\sin\theta(t) \tag{22}$$

where $R_y$ represents a vertical reaction force appearing between the capsule and underlying surface (contact load) and $R_x$ describes a horizontal force resulting from the motion of the pendulum that acts on the capsule.

In the setup under consideration, the pendulum angle $\theta(t)$ is the controlled quantity. Provided that the control system effectively rejects perturbations of $\theta(t)$, which result from the acceleration of the capsule and other factors, equation (20) can be omitted – the function $\theta(t)$ becomes an a priori known input signal.

$F_x$ is the friction force present in the system. Its definition (23) assumes the classical dry friction (Coulomb) model.

$$F_x(t) = \begin{cases} \mu R_y(t)\mathrm{sgn}[\dot{x}(t)] \leftrightarrow \dot{x}(t) \neq 0 \\ \mu R_y(t)\mathrm{sgn}[R_x(t)] \leftrightarrow \dot{x}(t) = 0 \wedge |R_x(t)| \geq \mu R_y(t) \\ R_x(t) \leftrightarrow \dot{x}(t) = 0 \wedge |R_x(t)| < \mu R_y(t) \end{cases} \tag{23}$$

Dimensionless parameters are introduced.

$$\Omega = \sqrt{\frac{g}{l}}, \tau = \Omega t, \gamma = \frac{M}{m}, z = \frac{x}{l},$$

$$f_z = \frac{F_x}{m\Omega^2 l}, r_z = \frac{R_x}{m\Omega^2 l}, r_y = \frac{R_y}{m\Omega^2 l} \tag{24}$$

Dependences between derivatives with respect to the dimensional and dimensionless time ($t$ and $\tau$ respectively) are presented in expression (25).

$$\dot{x} = \frac{dx}{dt} = \frac{dx}{d\tau}\frac{d\tau}{dt} = \Omega\frac{dx}{d\tau} = \Omega x'$$

$$\ddot{x} = \frac{d^2x}{dt^2} = \frac{d}{dt}\left(\frac{dx}{dt}\right) = \frac{d}{d\tau}\left(\Omega\frac{dx}{d\tau}\right)\frac{d\tau}{dt} = \Omega^2\frac{d^2x}{d\tau^2} = \Omega^2 x'' \tag{25}$$

Substitution of the relationships (23-25) into the equations (20-22) yields motion equations of the system under consideration (Fig. 1) in the dimensionless form.

$$z''(\tau)(\gamma+1) - \theta''(\tau)\cos\theta(\tau) = -\theta'^2(\tau)\sin\theta(\tau) - f_z(\tau) \tag{26}$$

$$r_y(\tau) = (\gamma+1) - \theta''(\tau)\sin\theta(\tau) - \theta'^2(\tau)\cos\theta(\tau) \tag{27}$$

$$r_z(\tau) = \theta''(\tau)\cos\theta(\tau) - \theta'^2(\tau)\sin\theta(\tau) \tag{28}$$

$$f_z(\tau) = \begin{cases} \mu r_y(\tau)\mathrm{sgn}[z'(\tau)] \leftrightarrow z'(\tau) \neq 0 \\ \mu r_y(\tau)\mathrm{sgn}[r_z(\tau)] \leftrightarrow z'(\tau) = 0 \wedge |r_z(\tau)| \geq \mu r_y(\tau) \\ r_z(\tau) \leftrightarrow |r_z(\tau)| < \mu r_y(\tau) \end{cases} \tag{29}$$

It is essential to note that in the scope of this paper, the goal is to optimize the function $\theta(\tau)$. Therefore, in each execution of the system (26-29), $\theta(\tau)$ and its derivatives with respect to the dimensionless time, are assumed to be known.

### 3.2 Designing of the pendulum capsule drive

Fig. 5 depicts a 3D CAD-rendered model of the pendulum capsule drive. The main aim of the research is to validate the efficiency of the proposed Fourier series-based algorithm in controlling the aforementioned system. Thus, the simplified non-encapsulated shape of the device's design is proposed.

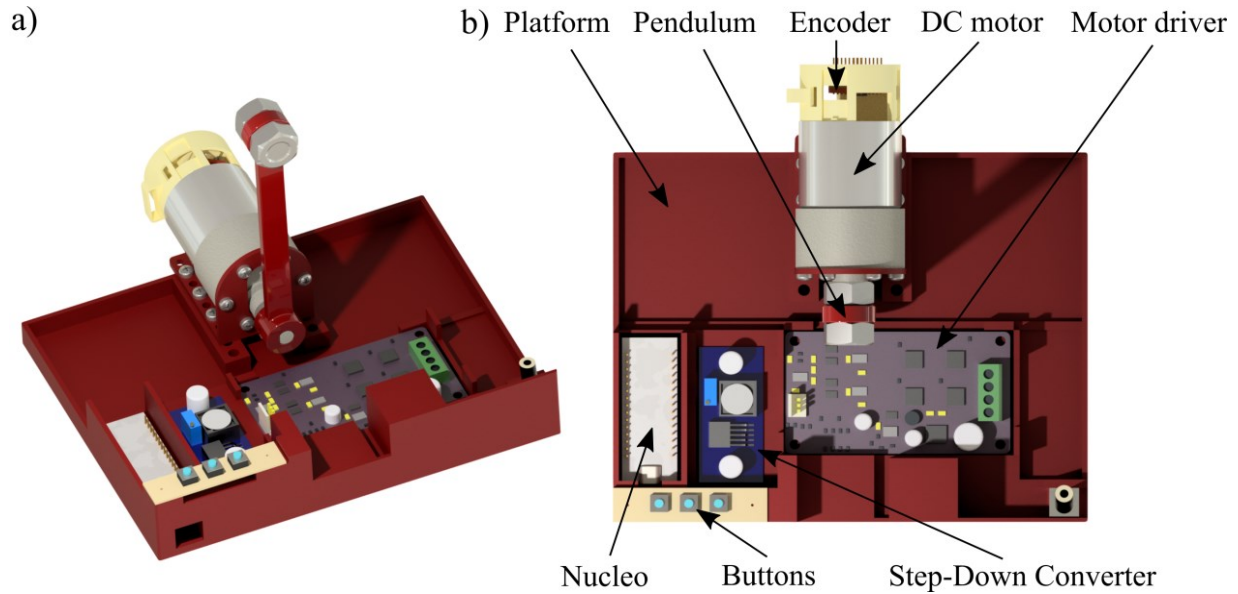

*Fig. 5 Pendulum capsule drive model: **a** general view; **b** top view*

The pendulum capsule drive's base platform measures 130 x 160 x 10 mm in size. The executive part of the system is a pendulum (length: 100 mm, mass: 40 g), attached to a DC motor POLOLU 4692 (max torque: 1.5 Nm, max rotational speed: 330 rpm, gear ratio: 30:1, operating voltage: 24V). The motor is powered via a CYTRON MD10C DC driver utilizing the Pulse Width Modulation (PWM) signal capable of controlling the direction and the speed of the motor. The control of the pendulum's position during its motion, as well as the generated signal tracking, is provided by the magnetic rotary encoder AS5040 (resolution: 10-Bit 360°) located on the rear of the motor. Additionally, the system employs a DC/DC step-down LM2596 converter to supply power to the STM32 Nucleo-32 board. The microcontroller handles all aspects of the device, while communication with the system, including behavioral assessment of the elements and data transmission, is facilitated by the ST-LINK acting as a UART-USB converter. On the upper section of the platform, three tact switches (buttons) designated for system calibration, program selection and launch, and device shutdown are positioned.

### 3.3 Control electronics

The subsequent pivotal step in ensuring the seamless operation of the system involved programming all electronic components and determining the device's control strategy. Consequently, this section elucidates the essential software aspects needed to achieve satisfactory results.

#### 3.3.1 PID control

The aim of the control strategy for the device is to follow the reference signal (control function calculated in the optimization process) represented by the positions of the pendulum in a specified interval of time. These values are sent to the DC motor by the driver through the Pulse Width Modulation (PWM) pin, determining the amount of energy required to achieve the target position by the pendulum. The available PWM range was adjusted as a trade-off between the PWM signal frequency and the resolution affecting its computations. Thus, to protect the human from the audible displeasing squeak of magnetic components and limit the switching losses of the power transitions in the motor driver, the frequency of 20 kHz was set, achieving a satisfactory PWM resolution of 3200 steps.

The closed-loop feedback control makes the realization of the target signal tracking feasible. The feedback is realized by the encoder, indicating the current position of the pendulum. The difference between the target pendulum position $\theta_z(t)$ (reference signal) and the position measured by the encoder $\theta(t)$ (output signal) determines the control error $e(t) = \theta_z(t) - \theta(t)$ for the feedback implementation. The constant sampling period for the calculations of integral and derivative parts of PID is maintained by the microcontroller's timer. The frequency of the PID loop is set at 100 Hz, providing the encoder with a satisfactory number of generated impulses (up to 845 for one pass of the PID loop) at the highest speed of the pendulum for valid and smooth regulation. To detect the motion, one input of the encoder is set as an interrupt, whereas the other is used for its direction identification (clockwise/counterclockwise).

Since the system is highly nonlinear, additional expressions compensating the friction and gravity interaction need to be considered when computing output values $u$. The scheme of the closed-loop feedback control is presented in **Fig. 6**.

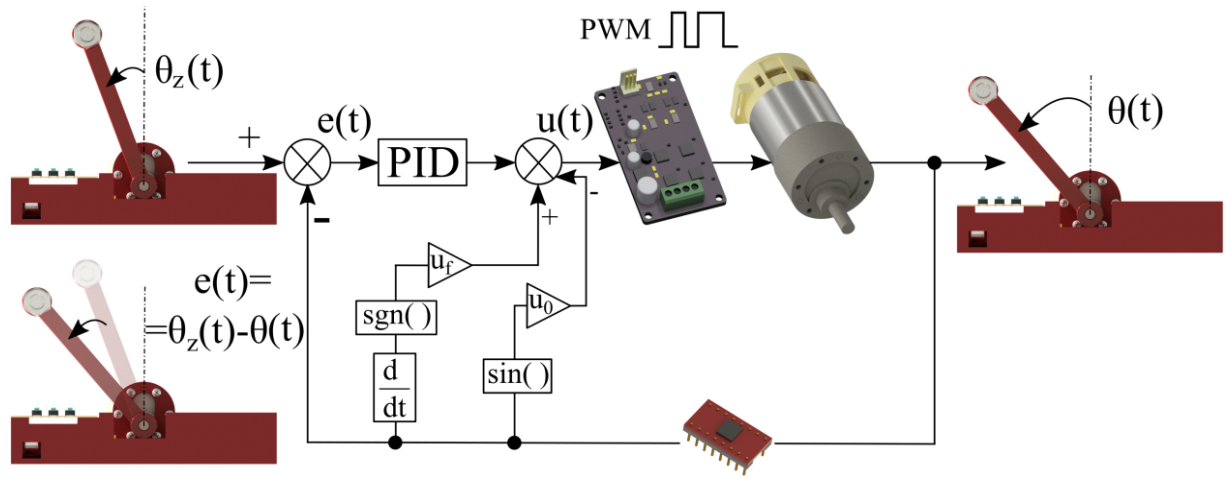

**Fig. 6** Closed-loop feedback control for the pendulum capsule drive. $\theta_z(t)$ –target pendulum position, $\theta(t)$ –measured position, $e(t)$ – control error, $u(t)$ – control signal, $u_f$ – friction compensation coefficient, $u_0$ – gravity compensation coefficient

The *friction compensation,* described by the expression:

$$u_f sgn(\dot{\theta}(t)) \tag{30}$$

accounts for a torque that is necessary to overcome friction occurring during the motion.

The friction compensation coefficient, denoted by $u_f$, is measured empirically. By gradually increasing the PWM duty cycle, the first value at which the pendulum swings out of its vertical position without any oscillations is identified.

The *gravity compensation* expressed as follows:

$$u_0 \sin(\theta(t)) \tag{31}$$

counteracts the nonlinearity induced by the gravity torque that affects motion of the pendulum.

The gravity compensation coefficient, denoted by $u_0$, is also measured empirically. By progressively increasing the PWM duty cycle, the first value that causes the pendulum to swing from the horizontal position without oscillations is registered. It should be noted that during the pendulum's swing towards its vertical position, frictional forces must also be overcome. Therefore, the previously determined coefficient $u_f$ needs to be subtracted. During the empirical measurement of the friction and gravity compensation coefficients, the PID controller must remain inactive.

### 3.3.2 User interface

Developing the user interface simplifies the requirement for manually uploading and launching the program for diverse signals and precise calibration procedures each time. For this, three tact switches are positioned on the upper part of the platform, integrating their functions with the LED3 available on the microcontroller. To provide an immediate reaction to the buttons' triggering, they are set as external interruptions and programmatically secured from debouncing.

Triggering of tact switches is responsible for the calibration procedure, the generated signal selection, and the chosen signal's launching, respectively. During the execution of any mode, pressing any of the switches again executes the Panic Button, which shuts down the system, protecting it from unpredicted states and potential damage.

## 3.4 Experimental set-up

To obtain the essential physical parameters of the device, the frame of the pendulum and the platform were fabricated using 3D printing technology. Transparent Polyethylene Terephthalate Glycol (PET-G) material was chosen for its strength and durability. The top part of the pendulum consists of two M10 stainless steel nuts with 24 mm long bolts inside. The pendulum is attached to the DC motor and placed on the platform, along with all the remaining electronic components making up the system, as described in Section 3.3. Underneath the platform, two skids made of bookbinding cardboard measuring 160x16x2 mm are fixed. To enable the device's motion, two aluminum profiles with a maximum length of 1.2 m were utilized to construct a running track. The view of the prototype on the track is shown in **Fig. 7**. The outer edge of the track is secured by aluminum bands, limiting the friction during potential rotations of the device. The gap between the bands was so set to ensure the platform's parallel guidance along the path. To prevent the device from leaping up due to higher generated inertia forces than expected, an additional external load of 50 g is included in the system.

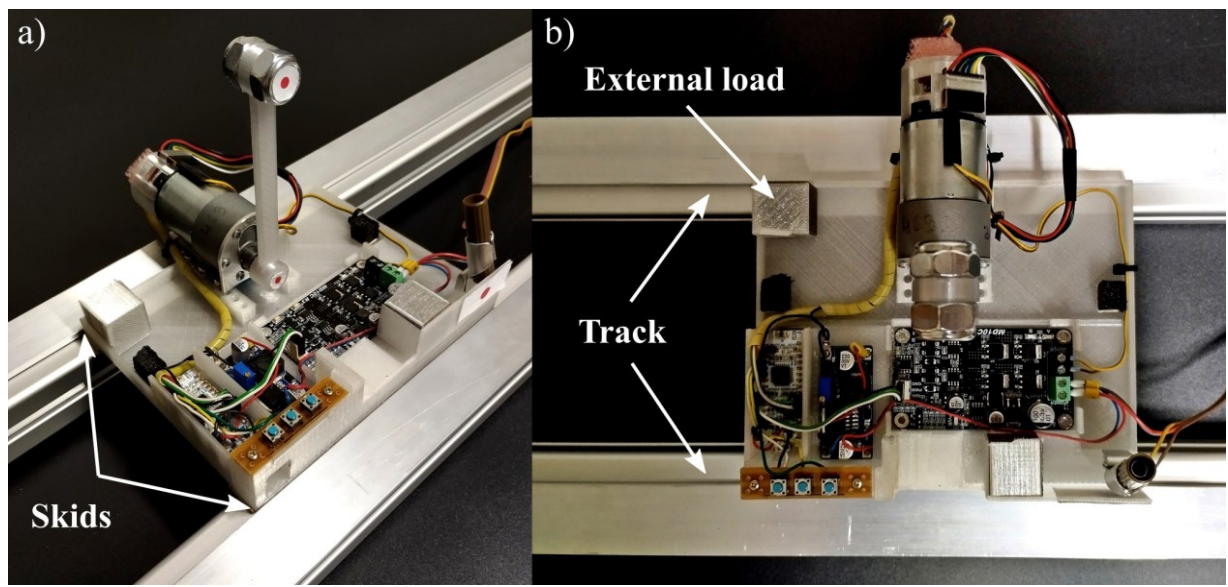

*Fig. 7 Prototype demonstration on the track: **a** general view; **b** top view*

Fig. 8 illustrates the experimental setup for the locomotion performance validation and the proposed method of controlling the device.

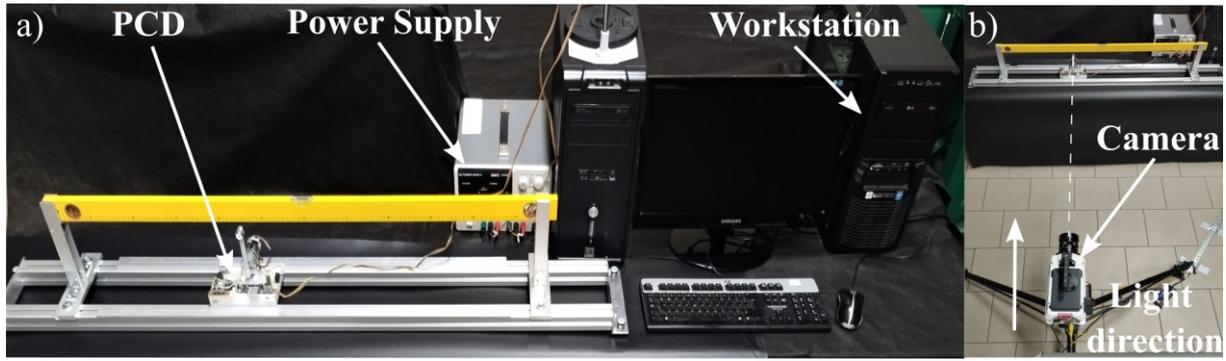

**Fig. 8** Experimental setup for pendulum capsule drive (PCD): **a** general view; **b** camera view

The experimental setup comprises the pendulum capsule drive positioned on the track, along with a measuring scale serving as a reference for the system's locomotion performance (see Fig. 8a). The device is operated from the power supply. The motion of the pendulum and the platform are captured with the use of a high-speed Phantom v711 camera with 900fps at full 1280x800 pixels monochromatic resolution, providing a smooth representation for the motion tracking software (see Fig. 8b). The camera is positioned 1.3 m away from the track. A Carl Zeiss planar lens 50 mm/f2 aimed at the mid-position of the reference distance covered by the system based on known numerical findings. This setup allows the system to capture the entire motion while minimizing optical distortions. In addition, a continuous light source maintains a stable background, preventing sudden shifts in nearby image elements.

## 3.5  Control optimization and its constraints

In this paper, the numerical investigation's aim is to optimize the rotational motion profile of the pendulum so that the system's total covered distance is maximum in a fixed interval of time. In this case, the control signal $u(\tau)$ corresponds to the state variable $\theta'(\tau)$, i.e., the controlled quantity is the rotational speed of the pendulum. Thus, the selection of the control function $u(\tau)$ directly affects not only the rotational speed of the pendulum $\theta'(\tau)$, but its angle $\theta(\tau)$ and angular acceleration $\theta''(\tau)$ as well. Consequently, choice of the appropriate control must be consistent with constraints imposed on all three quantities: $\theta(\tau), \theta'(\tau), \theta''(\tau)$. Restrictions regarding the angle $\theta(\tau)$ involve geometry of the laboratory stand and protection against pendulum hitting the body of the device. For these reasons, it is asserted that $\theta(\tau) \in \left[-\frac{\pi}{3}, \frac{\pi}{3}\right]$. Furthermore, since the pendulum is installed directly on the motor's shaft, the rotational speed of the pendulum must not exceed the maximum value for the motor, specified in its datasheet. Introducing the dimensionless quantities, the pendulum's speed must remain within the following interval: $\theta'(\tau) \in [-3.4, 3.4]$.

Apart from the kinematic properties of the pendulum's trajectory, the limited ability of the motor to generate a torque has to be taken into account, as formula (32) indicates:

$$\left|M_{max} - \zeta \dot{\theta}(t)\right| \geq \left|ml^2 \ddot{\theta}(t) + mgl \sin \theta(t)\right| \qquad (32)$$

where: $\zeta = \frac{M_{max}}{\omega_{max}}$, $M_{max}$ – maximum (stall) torque of the DC motor, $\omega_{max}$ – maximum rotational speed of the DC motor's shaft. The left-hand side of the inequality depicts the maximum output torque of the motor when the speed of its shaft equals $\dot{\theta}(t)$. The right-hand side represents the loads resulting from the presence of the pendulum, describing the influences of inertia and gravity, respectively. Division by $mgl$ and the introduction of the non-dimensional quantities $u_{max} = \frac{M_{max}}{mgl}$, $\kappa = \frac{\zeta \omega}{mgl}$, $\Omega = \sqrt{g/l}$, and the dimensionless time $\tau = \Omega t$ transform expression (32) into (33):

$$|u_{max} - \kappa\theta'(\tau)| \geq |\theta''(\tau) + \sin\theta(\tau)| \tag{33}$$

where $\theta' = \frac{d\theta}{d\tau}$. This condition must be fulfilled at all times; otherwise, the motor is not able to generate a sufficient torque to produce the required motion. Using the motor's datasheet and parameters of the device, the following values were obtained: $u_{max} = 25, \kappa = 10.85$. To ensure that the required torque can be physically achieved, a 30% safety margin is introduced for condition (33).

Since the centrifugal force acting on the pendulum may cause the platform to leap, affecting locomotion performance and signal tracking accuracy, a numerical limitation is provided to protect against this behavior. The leaps can occur when the centrifugal force $m\left(\dot{\theta}(t)\right)^2 l$ matches or exceeds the total weight of the pendulum and the platform: $(1 + \gamma)mg$. Therefore, at any moment, condition (34) should be fulfilled.

$$(1 + \gamma)mg \geq m\left(\dot{\theta}(t)\right)^2 l \tag{34}$$

Here, the worst-case scenario, in which the centrifugal force acts vertically upwards, is taken into account.[4] Dividing both sides of this inequality by $mgl$, introducing the parameter $\Omega = \sqrt{g/l}$ and the dimensionless time $\tau = \Omega t$ transforms the condition (34) into the form (35).

$$(1 + \gamma) \geq \left(\theta'(\tau)\right)^2 \rightarrow \theta'(\tau) \in \left[-\sqrt{1+\gamma}, \sqrt{1+\gamma}\right] \tag{35}$$

For the optimization procedure and then the simulation, the following values of the pendulum capsule drive prototype were identified: $\mu = 0.17$ and $\gamma = 14.5$. The static friction coefficient $\mu$ was measured empirically by gradually increasing the tilt angle until the device started to slide on the track and assumed as an approximated value of kinetic one. Note that for the adopted value of $\gamma$, condition (35) is less restrictive than the previously assumed limitation resulting from the motor's speed limit: $\theta'(\tau) \in [-3.4, 3.4]$.

The rotational motion profile parameters are optimized using the Differential Evolution (DE) algorithm to ensure that the pendulum capsule drive covers the greatest possible distance in a fixed interval of time. Thus, the cost function to be minimized is represented by formula (36):

$$J = -|z(\tau_f) - z(\tau_0)| \tag{36}$$

which involves minimizing the negative distance $z$.

The numerical computations are initialized optimizing the control function $u(\tau)$ for 3 harmonics with a period equal to $T_n = \frac{2\pi}{\omega_n}$, where $\omega_n = 1$. It should be noted that in this research, the value of the parameter $\omega$ is fixed in each iteration and not optimized. Once the result is achieved, it is used as an initial guess for the next iteration of optimization, represented by 6 harmonics with a doubled period $T_{n+1} = \frac{2\pi}{\omega_{n+1}}$, where $\omega_{n+1} = 0.5$, and 12 harmonics with a period $T_{n+2} = \frac{2\pi}{\omega_{n+2}}$, where $\omega_{n+2} = 0.25$.

### 3.6 Measurement technique

The image-based measurement method was employed to analyze the system's behavior, offering a viable alternative to traditional methods that require the installation of additional elements (e.g., measurement devices and wires) that unfavorably impact the system dynamics [55].

---

[4] Obviously, only the normal component of the acceleration of the mass center of the pendulum is considered here. This simplification results from the fact that this component dominates over tangential acceleration in all tested trajectories in the system. Therefore, the tangential component can be omitted in these considerations.

To track the motion of the pendulum and the platform, two circular markers are placed at characteristic points on the pendulum (see Fig. 9). The background is kept dark to prevent the effects of light reflections on recordings. The first marker ($O_1$) serves as a fixed reference point on the motor shaft, correlating with the pendulum's motion and indicating the platform's position. Simultaneously, the second marker ($O_2$) is placed on the movable upper section of the pendulum. These markers form the basis for the recognition stage in the open-source Kinovea tracking software. Experimental data is recorded as a movie consisting of a long sequence of consecutive images that track the position of each marker. The positions are digitally converted into Cartesian coordinates (x, y), measured in pixels and recalculated to millimeters. A detailed description of the Kinovea tracking procedure can be found in [55].

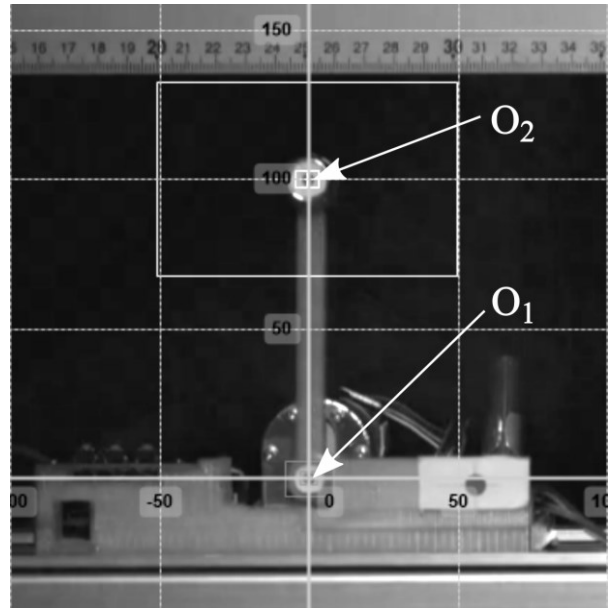

*Fig. 9 Setup of the markers for pendulum capsule drive motion tracking (snapshot from the Kinovea software)*

## 4   Results and discussion

The main goal of this work is to validate the robustness and efficiency of the pendulum capsule drive under Fourier series-based control, focusing on two parallel factors: forcing signal tracking accuracy and distance covered by the device in a fixed interval of time. The first step towards achieving this objective was to optimize the control function $u(\tau)$ over a short, fixed interval of dimensionless time for the chosen numbers of harmonics $k$ ($k = 3$, $T_1 = 2\pi$; $k = 6$, $T_2 = 4\pi$; $k = 12$, $T_3 = 8\pi$). The parameters of control functions $u(\tau)$ determined in this process were subsequently used in in simulations of motion of the pendulum capsule drive over 15.1 seconds. This value arises from the fact that within 15.1 seconds, a full number of periods of the control function can be contained (24 periods for $k = 3$, 12 for $k = 6$, 6 for $k = 12$). The achieved numerical values for distance and the corresponding recalculated speeds are presented in Table 1.

*Table 1* The comparison of the distance covered by the pendulum capsule drive in 15.1 s in the numerical and experimental stages, along with the speed and relative change for different numbers of harmonics k

| $k$ | Numerical Distance [cm] | Experimental Distance [cm] | Numerical Speed [cm/s] | Experimental Speed [cm/s] | Relative Difference δ [%] |
|---|---|---|---|---|---|
| | | | Time: 15.1 s | | |
| 3 | 37.8 | 37.4 ± 0.3 | 2.50 | 2.48 | 1.0 |
| 6 | 39.6 | 39.0 ± 0.4 | 2.62 | 2.58 | 1.7 |
| 12 | 39.6 | 39.0 ± 0.4 | 2.62 | 2.58 | 1.7 |

Analyzing the results, firstly, we observe that the analyzed control optimization algorithm based on Fourier series enabled the device to cover a distance of nearly 38 cm within the specified time of 15.1 seconds without relying on any data regarding the dynamics of the controlled system (except for optimization constraints). Secondly, according to the adopted assumptions regarding the greedy variant of the analyzed algorithm, the achieved result is at least non-decreasing with increasing number of harmonics $k$ in the optimized control signal. Thirdly, changing the number of harmonics from $k = 3$ to $k = 6$ results in noticeable improvement in control quality. Fourthly, we notice that increasing $k$ from 6 to 12 does not lead to any improvement in control quality. This is because the optimization algorithm for $k = 12$ remains at the initial conditions (6 harmonics from the previous iteration, followed by 6 zeros), resulting in an identical control signal for $k = 6$ and $k = 12$. Such phenomenon may occur due to the excessive dimensionality of the optimized parameter space or the optimization algorithm reaching its limits for the system under consideration at $k = 6$. Nevertheless, this result underscores the guarantee that results do not worsen with increasing harmonics. Of course, this fact also implies that further increasing $k$ in simulations is futile.

Several locomotion trials of the pendulum capsule drive were recorded under generated signals for $k = 3$ and $k = 6$ on the test stand, and subjected to comprehensive data analysis using Kinovea software, OriginPro, and Python scripts. Locomotion performance was assessed based on tracking point $O_1$. Experimental distances are consistent with numerical findings (see Tab. 1) and reproducible, as indicated by the low standard deviation (SD). Furthermore, it has been demonstrated that the assumption that higher numbers of harmonics improve locomotion performance is valid. The average speed of the device was determined based on these distances, and the relative difference between numerical and experimental findings was calculated using formula (37):

$$\delta = \frac{|v_{num} - v_{exp}|}{v_{num}} \cdot 100\% \tag{37}$$

where $v_{num}$ and $v_{exp}$ are the average speeds of the drive in the numerical simulation and the experiment, respectively. The relative errors in the experiment are only 1% and 1.7% compared to theoretical calculations. The comparison of locomotion performance, incorporating both numerical and experimental results for $k = 3$ (see Fig. 10a) and for $k = 6$ (see Fig. 11a), demonstrates the consistency of the findings, with minor observable differences. They may result from the assumption of the constant friction coefficient, micromechanical backlashes in the pendulum fixing, and mechanical gearbox backlashes. Both trajectories share a similar pattern of motion in which the stick-slip phenomenon occurs (due to the dry friction). In the stick phase, the velocity and the acceleration of the capsule are equal to zero, meaning that the frictional forces have not yet been disrupted or the inertia forces generated by the pendulum are too low. Hence, the pendulum capsule drive remains stationary for a very brief period of time. In the slip phase, the frictional forces have been disrupted and the inertia forces are sufficiently strong to cause the device's motion. The forward or backward drift of the device is determined by the pendulum motion, with the horizontal force playing a significant role. The motion stages of the pendulum capsule drive are presented in Fig. 12 and Fig. 13 for $k = 3$ and $k = 6$, respectively.

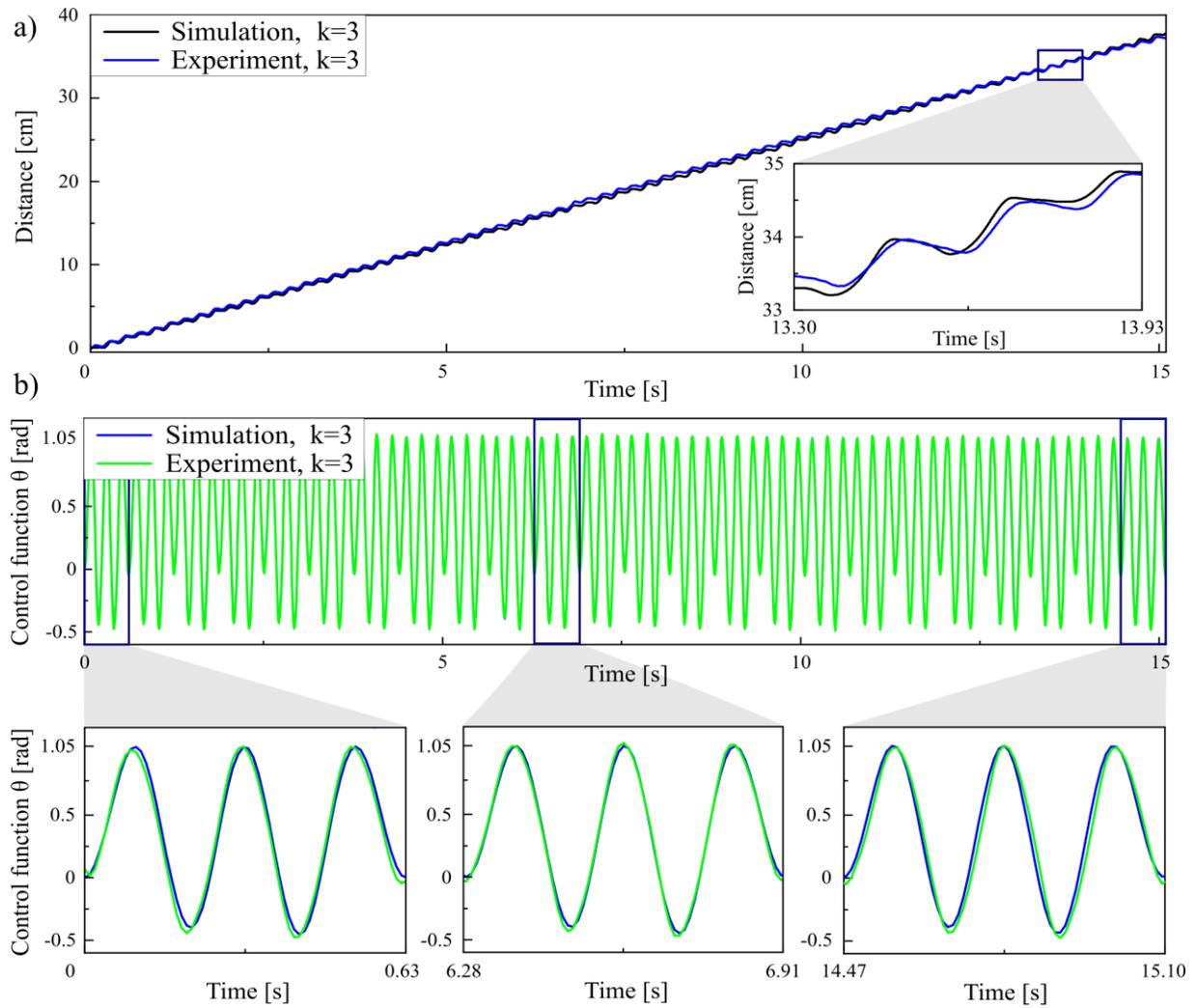

**Fig. 10** *Comparison of numerical and experimental results for 3 harmonics* ($k$)*: **a** locomotion performance; **b** control function (signal)*

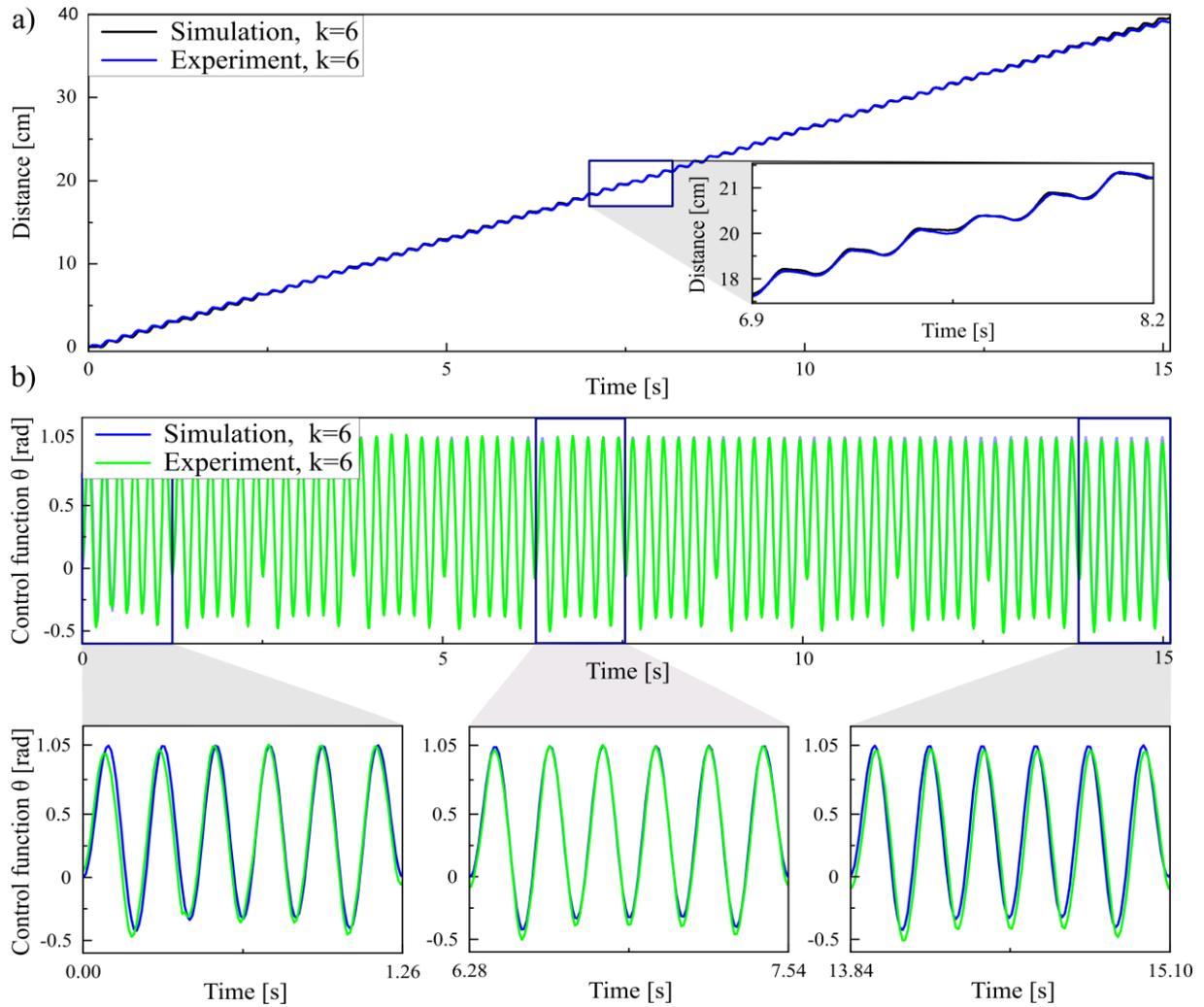

**Fig. 11** Comparison of numerical and experimental results for 6 harmonics ($k$): **a** locomotion performance; **b** control function (signal).

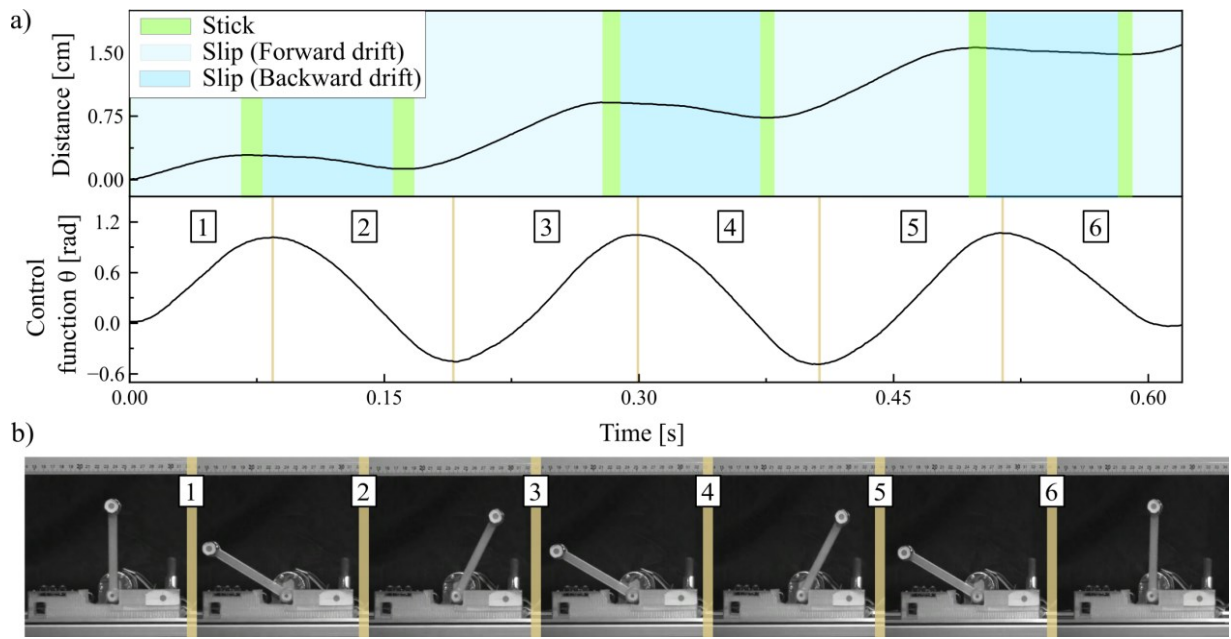

**Fig. 12** Pendulum capsule drive motion pattern for 3 harmonics ($k$): **a** stick-slip phases; **b** snapshots of the device motion

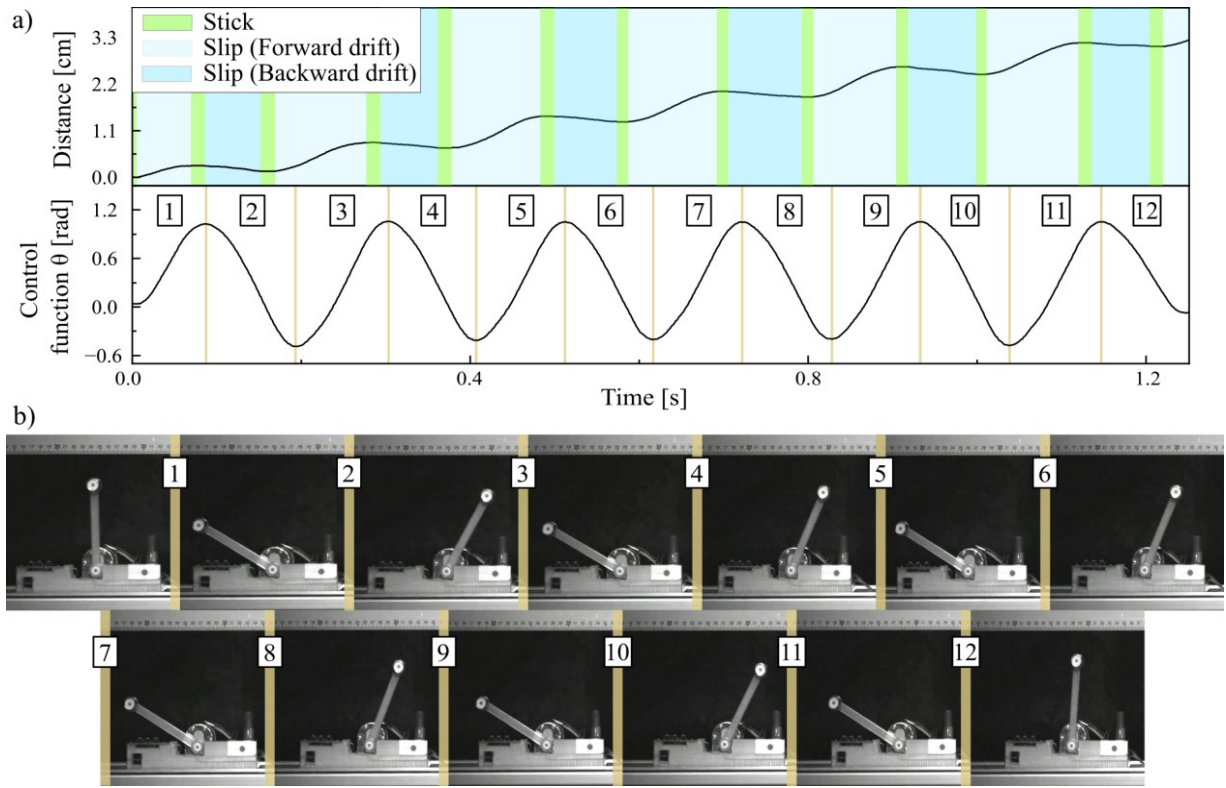

*Fig. 13* Pendulum capsule drive motion pattern for 6 harmonics ($k$): *a* stick-slip phases; *b* snapshots of the device motion

The angular displacement signal $\theta(t)$ was acquired using data from two markers, according to formula (38):

$$\theta = \operatorname{atan}\left(\frac{x(O_1) - x(O_2)}{y(O_2) - y(O_1)}\right) \quad (38)$$

where:

$O_1$ – tracking point of the bottom marker,

$O_2$ – tracking point of the upper marker.

This stage of the research involved a comprehensive assessment of signal tracking accuracy (averaging signals from recorded trials), calculated using the root mean square error (RMSE) metric— a powerful tool for determining differences between numerical and experimental values of the control function. Initially, the signal was analyzed over complete periods, with RMSE calculated for each period. Subsequently, the analysis extended over the entire 15.1 seconds to assess long-term signal tracking accuracy. Detailed RMSE values for both cases are provided in Table 2, while Fig. 10b and Fig. 11b visually represent these findings, demonstrating the quality of signal reproduction.

*Table 2* The values of root mean square error (RMSE) for 3 and 6 harmonics considering the entire signal and each period

| $k$ | $RMSE$ from 15.1 s [rad] | average $RMSE$ from the period [rad] |
| --- | --- | --- |
| 3 | 0.045 | 0.041 ± 0.017 |
| 6 | 0.083 | 0.056 ± 0.014 |

The presented results of the accuracy of the signal tracking are satisfying, which is reflected in the low values of the $RMSE$. Moreover, the quality for the entire signal is only slightly lower than those mentioned above (average period: $RMSE < 0.06$, entire signal: $RMSE < 0.09$) providing the system's robustness for potential perturbations during the motion on a long-term scale.

## 5 Summary and conclusions

This paper investigates a novel, Fourier series-based method of control optimization in a pendulum capsule drive. This system, renowned for its complex dynamics and discontinuities arising from dry friction, requires a particularly efficient and precise control method to enhance locomotion performance. The robustness and efficacy of the proposed control optimization method are thoroughly evaluated through numerical and experimental studies. The paper first outlines the method's original mathematical formulation, and its modified variant, which implements a greedy approach. The theoretical background is followed by a detailed exploration of the pendulum capsule drive, including mathematical modeling, designing, development of control electronics, control strategy, and the user interface. The experimental setup, parameter identification, and adaptation of measurement technique are then meticulously detailed, followed by a presentation of the results.

The introduced control optimization method proposes an easy-to-deploy and efficient way of estimating open-loop optimal control incorporating constraints that are often encountered in real-world scenarios, such as a maximum range of available torque and speed of the DC motor. The algorithm's flexibility is provided by specifying minimal information about the object, a set of admissible controls, constant control ranges, and a unique performance for measuring any admissible control. Through the simplification of the optimized function to nonlinear programming (parameter optimization), the Differential Evolution (DE) algorithm can be used to find an optimized solution. However, solving problems with DE induces challenges while the number of parameters in the optimized control function increases, resulting in the growth of the number of iterations essential to finding an acceptable problem-solving outcome. Consequently, its greedy variant is proposed to counteract the vulnerability of the original Fourier series-based method. In this approach, the problem is split into smaller parts (short time intervals) for which the function is optimized, and then a starting point is set for optimizing over a more extended period until the ultimate objective across the entire domain is reached. Consequently, the previously found solution of control optimization is considered in subsequent iteration, ensuring a non-decreasing quality of solutions.

The algorithm can handle the discontinuous system that has been proven in the numerical investigation, demonstrating an effective way of controlling the locomotion of the pendulum capsule drive. Moreover, the increase of the non-decreasing quality of the solution (locomotion performance) as the number of harmonics $k$ in the control function is confirmed while optimizing the rotational motion profile. The aforementioned findings were validated on the experimental setup, where the locomotion performance and signal tracking accuracy were assessed in parallel. The experimental results are consistent with the numerical findings ($< 2\%$ loss), showcasing the algorithm's practical applicability and repeatability. For $k = 3$, the pendulum capsule drive achieved an average speed of 2.48 cm/s, while for $k = 6$, this increased to 2.58 cm/s, confirming the theoretical assumption of reaching longer distances with the inclusion of a higher number of harmonics in the optimized control function. The locomotion performance results are supported by satisfying signal-tracking accuracy, with an average root mean square error $RMSE < 0.06$ for each period and the high resolution of the Fourier signal reproduction over a long-term scale (15.1 s) with $RMSE < 0.09$, indicating resistance to potential perturbations during motion. Minor differences in experimental results may arise from the assumption of a constant friction coefficient, micromechanical backlashes in the pendulum fixing after several runs, and gearbox backlashes.

This paper confirms that the Fourier series-based method of control optimization can be effectively used in nonlinear and discontinuous systems, which is supported by numerical and experimental findings. The results indicate a potential pathway for future research in control optimization of capsule robots operating in complex environments, such as endoscopic procedures, rescue operations in confined spaces, or pipeline inspections.

## Author contributions

**Sandra Zarychta:** Conceptualization, Data curation, Formal analysis, Funding acquisition, Investigation, Methodology, Project administration, Resources, Software, Supervision, Validation, Visualization, Writing – original draft, Writing – review & editing.

**Marek Balcerzak:** Conceptualization, Investigation, Methodology, Software, Writing – original draft, Writing – review & editing.

**Katarzyna Wojdalska:** Writing – original draft, Writing – review & editing.

**Rafał Dolny:** Resources, Software, Writing – review & editing.

**Jerzy Wojewoda:** Conceptualization, Data curation, Investigation, Methodology, Resources, Supervision, Validation, Writing – review & editing.

## Funding

This study has been supported by the National Science Centre, Poland, PRELUDIUM Programme (Project No. 2020/37/N/ST8/03448).

## Acknowledgments

This paper has been completed while the first author was the Doctoral Candidate in the Interdisciplinary Doctoral School at the Lodz University of Technology, Poland.

The authors thank Wiesław Prus for setting up, providing equipment and technical assistance during the experiments.

## Availability of data and materials

The datasets and materials generated and analyzed during the current study are available from the corresponding author upon reasonable request.

## Supplementary materials

Supplementary material associated with this article can be found in the online version.

## Declarations

## Conflict of interest

The authors declare no conflict of interest in preparing this article.

## Ethical approval

Not applicable.